\documentclass[10pt,a4paper,reqno]{amsart}

\usepackage[utf8]{inputenc}
\usepackage[T1]{fontenc}
\usepackage[english]{babel}

\usepackage{amsmath}
\usepackage{amsfonts}
\usepackage{amssymb}
\usepackage{amsthm}
\usepackage{ucs}
\usepackage{graphicx}
\usepackage{xcolor}
\usepackage{dsfont}
\usepackage{tikz}
\usepackage{wasysym}
\usepackage{physics}
\usepackage{mathtools}
\usepackage{xifthen}
\usepackage{todonotes}
\usepackage[inline]{enumitem}
\usepackage{stmaryrd}
\usepackage{constants}
\usepackage[ruled]{algorithm2e}
\usepackage{tikz,pgfplots}
\usepackage{chngcntr}
\usepackage{hyperref}
\usepackage[font={small}]{caption}

\counterwithin{figure}{section}

\graphicspath{{Pics/}}

\usepackage{constants}

\newconstantfamily{ctrlcst}{
	symbol=\kappa
}
\newconstantfamily{csts}{
	symbol=C
}
\renewconstantfamily{normal}{
	symbol=c
}

\makeatletter
\newcommand{\pushright}[1]{\ifmeasuring@#1\else\omit\hfill$\displaystyle#1$\fi\ignorespaces}
\newcommand{\pushleft}[1]{\ifmeasuring@#1\else\omit$\displaystyle#1$\hfill\fi\ignorespaces}
\makeatother

\newcommand{\Rd}{\mathbb{R}^d}
\newcommand{\Zd}{\mathbb{Z}^d}

\newcommand{\Zplus}{\mathbb{Z}_+}

\newcommand{\betac}{\beta_{\mathrm{\scriptscriptstyle c}}}
\renewcommand{\norm}[1]{\|#1\|}
\renewcommand{\abs}[1]{\lvert#1\rvert}
\newcommand{\normI}[1]{\left\|#1\right\|_{\scriptscriptstyle 1}}

\newcommand{\setof}[2]{\{#1\,:\,#2\}}
\newcommand{\bsetof}[2]{\bigl\{#1\,:\,#2\bigr\}}
\newcommand{\given}{\,|\,}

\newcommand{\icl}{\tau}
\newcommand{\iclb}{\icl_{\beta,0}}
\newcommand{\iclbh}{\icl_{\beta,h}}
\newcommand{\comp}{\mathrm{c}}
\newcommand{\eone}{\mathbf{e}_1}

\newcommand{\p}{\mathbb{P}}

\newcommand{\Ham}{\mathcal{H}}

\renewcommand{\nleftrightarrow}{\mathrel{\ooalign{$\leftrightarrow$\cr\hidewidth$/$\hidewidth}}}

\renewcommand{\emptyset}{\varnothing}

\newcommand{\calE}{\mathcal{E}}
\newcommand{\calF}{\mathcal{F}}
\newcommand{\calG}{\mathcal{G}}

\newcommand{\calN}{\mathcal{N}}

\newcommand{\calU}{\mathcal{U}}

\newcommand{\sfo}{\mathsf{o}}

\newcommand{\sfq}{\mathsf{q}}

\newcommand{\sfw}{\mathsf{w}}

\newcommand{\sfO}{\mathsf{O}}

\newcommand{\sfZ}{\mathsf{Z}}

\newcommand{\bbB}{\mathbb{B}}
\newcommand{\bbC}{\mathbb{C}}

\newcommand{\bbE}{\mathbb{E}}

\newcommand{\bbG}{\mathbb{G}}

\newcommand{\bbR}{\mathbb{R}}

\newcommand{\bbZ}{\mathbb{Z}}

\newcommand{\normII}[1]{\|#1\|_{\scriptscriptstyle 2}}
\newcommand{\supp}{\mathsf{S}}

\newcommand{\uvec}{\mathbf{u}}

\newcommand{\defby}{\stackrel{\mathrm{\scriptscriptstyle def}}{=}}
\newcommand{\smo}{\sfo}

\newcommand{\bfD}{\mathbf{D}}
\newcommand{\define}[1]{\textsf{#1}}
\newcommand{\difsymcent}[1]{\vcenter{\hbox{$#1\triangle$}}}
\newcommand{\difsym}{\mathbin{\mathchoice
        {\difsymcent\scriptstyle}{\difsymcent\scriptstyle}
        {\difsymcent\scriptscriptstyle}{\difsymcent\scriptscriptstyle}}}

\newcommand{\fcone}{\mathcal{Y}^\blacktriangleleft}
\newcommand{\bcone}{\mathcal{Y}^\blacktriangleright}

\newcommand{\RWP}{\mathsf{P}}
\newcommand{\RWE}{\mathsf{E}}

\newcommand{\isingLaw}[2]{\ifthenelse{\equal{#2}{}}{\langle#1\rangle_{\beta}}{\langle#1\rangle_{\beta,#2}}}

\newcommand{\evGraphSet}[2]{\ifthenelse{\equal{#2}{}}{\mathcal{E}_{#1}}{\mathcal{E}_{#1}(#2)}}

\newcommand{\evGraphLaw}[2]{\ifthenelse{\equal{#2}{}}{\mathfrak{P}^{#1}}{\mathfrak{P}^{#1}_{#2}}}

\newcommand{\HTPF}[2]{\ifthenelse{\equal{#2}{}}{\mathbf{Z}^{\rm\scriptscriptstyle HT}(#1)}{\mathbf{Z}^{\rm\scriptscriptstyle HT}_{#2}(#1)}}

\newcommand{\currentLaw}[2]{\ifthenelse{\equal{#2}{}}{\mathbf{P}^{#1}}{\mathbf{P}^{#1}_{#2}}}
\newcommand{\currentPF}[2]{\ifthenelse{\equal{#2}{}}{\mathbf{Z}^{\rm\scriptscriptstyle RC}(#1)}{\mathbf{Z}_{#2}^{\rm\scriptscriptstyle RC}(#1)}}

\newcommand{\RCMLaw}[2]{\ifthenelse{\equal{#2}{}}{\p\left(#1\right)}{\p_{#2}\left(#1\right)}}

\newcommand{\parityGraph}[2]{\ifthenelse{\equal{#2}{}}{\mathbf{G}_{#1}}{\mathbf{G}_{#1}(#2)}}

\newcommand{\bfn}{\mathbf{n}}

\newcommand{\bfm}{\mathbf{m}}
\newcommand{\weight}{\textnormal{w}}

\newcommand{\RCPF}{Z}
\newcommand{\evenPart}{\mathcal{E}}

\renewcommand{\emptyset}{\varnothing}

\theoremstyle{plain}
\newtheorem{theorem}{Theorem}[section]
\newtheorem{lemma}[theorem]{Lemma}

\newtheorem{remark}{Remark}[section]

\theoremstyle{definition}
\newtheorem{obs}{Observation}

\author{S\'{e}bastien Ott}
\address{Section de Mathématiques, Université de Genève, CH-1211 Genève, Switzerland}
\email{sebastien.ott@unige.ch}
\author{Yvan Velenik}
\address{Section de Mathématiques, Université de Genève, CH-1211 Genève, Switzerland}
\email{yvan.velenik@unige.ch}

\title{Asymptotics of correlations in the Ising model:\\a brief survey}

\begin{document}
	
\begin{abstract}
We present a brief survey of rigorous results on the asymptotic behavior of correlations between two local functions as the distance between their support diverges, concentrating on the ferromagnetic Ising model on \(\Zd\).
\end{abstract}

\maketitle

%
%

\section{Introduction}

\subsection{Density--density correlations in fluids}

The analysis of the correlations of a random field is a problem of central importance both in probability theory and in statistical physics.
Indeed, for a given random field on \(\Zd\) or \(\Rd\) satisfying suitable regularity conditions, the whole structure of the field is encoded in its \define{Ursell functions} (or \define{cumulants}). If \(X\) is a random field on \(\Rd\) of law \(\mu\) and \(f_1,\dots,f_n\) are \(n\) local functions of the field, the \(n^{\text{th}}\) Ursell function of \(f_1,\dots,f_n\) is defined by
\begin{equation*}
	U_n(f_1,\dots,f_n) \defby \frac{\partial^{n}}{\partial t_n\cdots\partial t_1} \log(\mu\big( e^{\sum_{i=1}^n t_i f_i}\big) ) \Big|_{t_1=\dots=t_n=0}.
\end{equation*}
Of particular importance are the first Ursell function (the expected value) and the second one (the covariance). For example, in the case of Gaussian fields, these two determine the field's law. For the case of covariances, investigations were initiated more than a century ago, in 1914, in very influential works by Ornstein and Zernike~\cite{Ornstein+Zernike-1914,Zernike-1916}. Let us briefly sketch the main ideas they introduced and one of the consequences they derived. We will present them in the language of the Ising model, despite the obvious anachronism (the Ising model was invented several years after these works).

Let us denote by \(G_\beta(0,x) \defby \mu_\beta^+(\sigma_0\sigma_x) - \mu_\beta^+(\sigma_0)\mu_\beta^+(\sigma_x)\) the truncated 2-point function of the model, that is, the covariance of the spins located at vertices \(0\) and \(x\) of the lattice \(\Zd\) under the infinite-volume Gibbs measure \(\mu_\beta^+\) (if needed, the reader can find precise definitions of these notions below). In order to understand \(G_\beta(0,x)\), Ornstein and Zernike made the following assumptions:
\begin{enumerate}[label=(\roman*)]
    \item\label{eq:OZ_asmptn:FinCorLen} there exists \(a>0\) such that \(G_\beta(0,x) \leq e^{-a\norm{x}}\) for all \(x\in\Zd\);
    \item\label{eq:OZ_asmptn:Renewal} there exists a function \(C_\beta:\Zd\times\Zd \to \bbR\), the \define{direct correlation function}, such that
        \begin{equation}\label{eq:OZ:NaiveRenewal}
            G_\beta(0,x) = C_\beta(0,x) + \sum_{y\in\Zd} C_\beta(0,y) G_\beta(y,x) 
        \end{equation}
        and, for some \(B,b\in(0,\infty)\),
    	\begin{equation}\label{eq:OZ_asmptn:SepMasses}
    		C(0,x) \leq Be^{-b\norm{x}} G(0,x)\qquad\forall x\in\Zd.
    	\end{equation}
\end{enumerate}
Assumption~\ref{eq:OZ_asmptn:FinCorLen} states that the correlation length is finite; in other words, it assumes that the system is not critical: \(\beta\neq\betac\). 
Assumption~\ref{eq:OZ_asmptn:Renewal} relies on the following intuitive picture: the correlation between the spins at \(0\) and \(x\) possesses two sources. On the one hand, there is a direct influence \(C_\beta(0,x)\) of the spin at \(0\) on the spin at \(x\); on the other hand, there is the direct influence \(C_\beta(0,y)\) of the spin at \(0\) on a third spin at \(y\) combined with the full influence \(G_\beta(y,x)\) of the spin at \(y\) on the spin at \(x\). This is the content of~\eqref{eq:OZ:NaiveRenewal}.
They then assumed that the direct influence of one spin on another should decay with a strictly larger exponential rate than the full influence. This is the \define{separation of masses} assumption appearing in~\eqref{eq:OZ_asmptn:SepMasses}.

Notice that Ornstein and Zernike did not provide an independent definition of the function \(C_\beta\). In their approach, it must then be considered as being defined by the relation~\eqref{eq:OZ:NaiveRenewal}. This is a weak point of the approach, which is usually complemented (in the physics literature) with a recipe to construct an approximate version of the function \(C_\beta\). There are several such ``closure relations'', of which  the Percus--Yevick approximation and the hypernetted-chain equation are two common examples~\cite{Hansen+McDonald-2013}.

Alternatively, it was proposed to construct an explicit function \(C_\beta\), based on suitable diagrammatic expansions, and then prove that the function thus constructed does indeed satisfy~\eqref{eq:OZ:NaiveRenewal}. This is the approach used by Abraham and Kunz in one of the first rigorous derivations of the Ornstein--Zernike (OZ) theory~\cite{Abraham+Kunz-1977}. The origins of graphical approaches to the analysis of correlation functions seems to be Symanzik's famous paper~\cite{Symanzik-1968}.

In any case, assuming that the above assumptions are met, Ornstein and Zernike were able to derive several important results. Of particular relevance for the present review is their derivation of the asymptotic behavior of the function \(G_\beta(0,x)\) for large \(\norm{x}\). A version of this derivation can be found in~\cite{Zernike-1916}, where it is shown that, roughly speaking,
\begin{equation}\label{eq:OZ:NaiveAsymptotics}
G_\beta(0,x) \sim \norm{x}^{-(d-1)/2}\, e^{-\tau_\beta \norm{x}} ,
\end{equation}
where \(\tau_\beta\in (0,\infty)\) is the inverse correlation length. The prediction that the correction to the exponential decay is a prefactor of order \(\norm{x}^{-(d-1)/2}\) turns out to be correct in a wide variety of situations and is now known as OZ decay. \eqref{eq:OZ:NaiveAsymptotics} can be derived from the assumptions above as follows. Introducing the generating functions \(\hat G_\beta(k) = \sum_x e^{k\cdot x} G_\beta(0,x)\) and \(\hat C_\beta(k) = \sum_x e^{k\cdot x} C_\beta(0,x)\), the renewal equation~\eqref{eq:OZ:NaiveRenewal} becomes (note that \(G_\beta\) is translation invariant)
\begin{equation}\label{eq:OZEqn}
    \hat G_\beta(k) = \frac{\hat C_\beta(k)}{1-\hat C_\beta(k)} .
\end{equation}
Using lattice symmetries and the above assumptions, one can then use the approximation \(\hat{C}_\beta(k) = c_0 + c_1 \norm{k}^2 + \sfO(\norm{k}^4)\) with \(1>c_0>0\) and \(c_1>0\). The decay~\eqref{eq:OZ:NaiveAsymptotics} then follows from a Tauberian theorem.

The theory developed by Ornstein and Zernike has played an important role in the development of the statistical theory of fluids and is still used in its modern incarnations. Its prediction~\eqref{eq:OZ:NaiveAsymptotics} was first verified in the 1960s through explicit computations in the planar Ising model, but was also falsified in some regimes for the same model, showing that the above assumptions do not necessarily hold (see, in particular, Section~\ref{sec:LTheuristic} below). Roughly a decade later, a stream of rigorous works based on a variety of approaches established rigorous versions of~\eqref{eq:OZ:NaiveAsymptotics} in various perturbative regimes (e.g., very high or very low temperature), as well as suitable corrections when OZ decay turns out not to hold. Only recently (roughly in the last 15 years), has it become possible to derive versions of such results in a nonperturbative manner. As will be explained below, these modern probabilistic approaches are still based on a suitable reinterpretation of the assumptions introduced in the original paper by Ornstein and Zernike, combined with various graphical representations of the correlation functions.

\subsection{The Ising model and basic quantities}
In this review, we focus on the ferromagnetic Ising model on \(\Zd\), \(d\geq 2\) (the issues we discuss are essentially trivial when \(d=1\)). For the sake of simplicity, we will assume nearest-neighbor interactions, although most results have been obtained for finite-range interactions, and refer the reader to the original papers for these generalizations.

\medskip

\define{Configurations} of the Ising model on \(\Zd\) are given by elements \(\omega=(\omega_i)_{i\in\Zd}\) of the set \(\Omega\defby\{\pm 1\}^{\Zd}\). We denote by \(\sigma_i\) the \define{spin} random variable \(\sigma_i(\omega)=\omega_i\). Given a subset \(\Delta\subset\Zd\), we denote by \(\calF_\Delta\) the \(\sigma\)-algebra generated by the random variables \(\sigma_i\) with \(i\in\Delta\) and write \(\calF\defby\calF_{\Zd}\). The set of nearest-neighbor edges of \(\Zd\) is denoted \(\calE\defby\bsetof{\{i,j\}\subset\Zd}{\normI{j-i}=1}\). Finally, given \(\Delta\subset\Zd\), we write \(\Delta^\comp\defby\Zd\setminus\Delta\).

\medskip

Given \(\beta>0\), \(h\in\bbR\) and a subset\footnote{\(A\Subset\Zd\) means that \(A\subset\Zd\) and the cardinality \(\abs{A}\) of \(A\) is finite.} \(\Lambda\Subset\Zd\), the \define{Hamiltonian} of the Ising model in \(\Lambda\) is the function
\[
\Ham_\Lambda \defby -\beta \sum_{\{i,j\}\in\calE_\Lambda} \sigma_i\sigma_j - h\sum_{i\in\Lambda} \sigma_i,
\]
where \(\calE_\Lambda\defby\bsetof{\{i,j\}\in\calE}{\{i,j\}\cap\Lambda\neq\emptyset}\).

The \define{finite-volume Gibbs measure} in \(\Lambda\), at parameters \(\beta\) and \(h\) and with \define{boundary condition} \(\eta\in\Omega\), is the probability measure \(\mu_{\Lambda;\beta,h}^\eta\) on \(\Omega_\Lambda^\eta\defby\setof{\omega\in\Omega}{\omega_i=\eta_i\;\forall i\not\in\Lambda}\) given by
\[
\mu_{\Lambda;\beta,h}^\eta(\omega) \defby \frac{e^{-\Ham_\Lambda(\omega)}}{\sfZ_\Lambda^\eta}, 
\]
where the normalization constant \(\sfZ_\Lambda^\eta\) is known as the \define{partition function}.
We denote by \(\mu_{\Lambda;\beta,h}^+\) and \(\mu_{\Lambda;\beta,h}^-\) the measures corresponding to the boundary conditions \(\eta^+\equiv 1\) and \(\eta^-\equiv -1\) respectively.

We denote by \(\calG(\beta,h)\) the set of all (infinite-volume) \define{Gibbs measures} corresponding to the parameters \(\beta\) and \(h\), that is, the set of probability measures \(\mu\) on \((\Omega,\calF)\) that satisfy the \define{DLR equations}:
\[
\mu(\,\cdot\given\calF_{\Lambda^\comp})(\omega) = \mu_{\Lambda;\beta,h}^\omega(\,\cdot\,)\qquad\text{for \(\mu\)-a.e. \(\omega\in\Omega\) and every \(\Lambda\Subset\Zd\).}
\]
It is well-known (see, for instance, \cite{Friedli+Velenik-2017}) that \(\calG(\beta,h)\) is a nonempty simplex for all \(\beta,h\); we denote by \(\mathrm{ext}\calG(\beta,h)\) the set of all \define{extremal measures}, that is, the set of extremal elements of \(\calG(\beta,h)\). The weak limits \(\mu_{\beta,h}^\pm\defby\lim_{\Lambda\to\Zd} \mu_{\Lambda;\beta,h}^\pm\) always belong to \(\mathrm{ext}\calG(\beta,h)\) (and are, in fact, the only elements of this set that are translation invariant). Moreover, non-uniqueness (that is, \(\abs{\calG(\beta,h)}>1\)) is equivalent to \(\mu_{\beta,h}^+(\sigma_0)\neq\mu_{\beta,h}^-(\sigma_0)\), which occurs if and only if \(h=0\) and \(\beta > \betac(d)\), where \(\betac(d)\in(0,\infty)\) for all \(d\geq 2\). In all other cases, there is a unique Gibbs measure (in particular, \(\mu_{\beta,h}^+=\mu_{\beta,h}^-\)).

The quantity \(\mu_{\beta,0}^+(\sigma_0)\) is known as the \define{spontaneous magnetization}; it satisfies \(\mu_{\beta,0}^+(\sigma_0)=-\mu_{\beta,0}^-(\sigma_0)\) and is positive for all \(\beta>\betac(d)\) and zero for all \(\beta\leq\betac(d)\).

\medskip

A function \(f:\Omega\to\bbR\) is \define{local} if there exists a finite set \(\Lambda\Subset\Zd\) such that \(f\) is \(\calF_\Lambda\)-measurable; the smallest such set if called the \define{support} of \(f\) and denoted \(\supp(f)\).
We will write \(\mu(f;g) \defby \mu(fg)-\mu(f)\mu(g)\) the covariance between two functions \(f\) and \(g\). When \(f=\sigma_i\) and \(g=\sigma_j\), the corresponding quantity is also known as the \define{truncated 2-point function}.

\medskip

It is a general property of extremal measures that local functions are asymptotically uncorrelated:
\[
\lim_{R\to\infty}\sup_{\substack{f,g \text{ local}\\d(\supp(f),\supp(g))>R}} \mu(f;g) = 0
\qquad\text{for all }\mu\in\mathrm{ext}\calG(\beta,h).
\]
As mentioned in the beginning of this section, an important problem in statistical mechanics (and probability theory) is to obtain more quantitative information about the speed of decay of these covariances. 

\section{Exponential decay}

The \define{inverse correlation length} is defined as the rate of exponential decay of the truncated 2-point function: given a unit vector \(\uvec\) in \(\Rd\),
\[
\iclbh(\uvec) \defby - \lim_{n\to\infty} \frac1n \log\mu_{\beta,h}^+(\sigma_0;\sigma_{[n\uvec]}),
\]
where \([x]\in\Zd\) denotes the component-wise integer part of \(x\in\Rd\). The existence of \(\iclbh\) can be established using a subadditive argument, based either on a suitable graphical representation of the truncated 2-point function~\cite{Graham-1982}, or on reflection positivity using a transfer operator (see~\cite{Seiler-1982}).

Of course, the existence of \(\iclbh\) does not imply that the truncated 2-point function actually decays exponentially fast, since \(\iclbh\) might in fact equal \(0\). It turns out, however, that the latter happens only at the critical point:
\[
\inf_{\uvec} \iclbh(\uvec) > 0 \text{ for all }(\beta,h)\neq (\betac(d),0) \qquad\text{ and }\qquad \icl_{\betac(d),0} = 0.
\]
This was shown, for the Ising model with finite-range interactions, in~\cite{Lebowitz+Penrose-1068} when \(h\neq 0\), in~\cite{Aizenman+Barsky+Fernandez-1987} when \(h=0\) and \(\beta<\betac(d)\), and in~\cite{Duminil-Copin+Goswami+Raoufi} when \(h=0\) and \(\beta>\betac(d)\). The result for \(h=0\) and \(\beta=\betac(d)\) follows, for the nearest-neighbor model, from a lower bound on the susceptibility proved in~\cite{Simon-1980} and the fact the spontaneous magnetization vanishes at \(\betac(d)\), proved in~\cite{Aizenman+Duminil-Copin+Sidoravicius-2015}.

\medskip

\begin{figure}[t]
    \includegraphics[width=10cm]{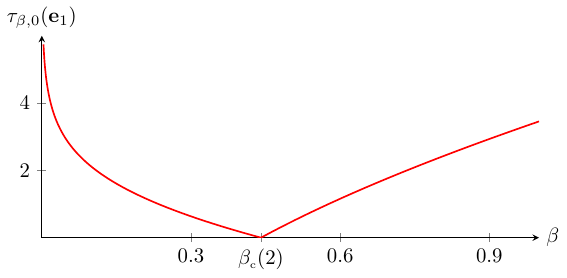}
    \caption{The inverse correlation length of the 2d Ising model at \(h=0\) as a function of the inverse temperature \(\beta\). Observe how it diverges as \(\beta\to 0\) ar \(\beta\to\infty\) and vanishes at \(\betac(2)=\frac12\mathrm{argsinh}(1)\cong 0.44\).}
    \label{fig:iclb2d-beta}
\end{figure}

\begin{figure}[t]
    \includegraphics[width=10cm]{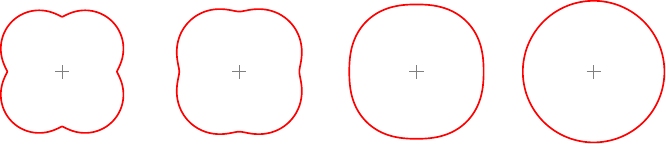}
    \caption{The inverse correlation length of the 2d Ising model as a function of the direction (for \(h=0\) and four increasing values of the temperature, all below the critical temperature). Only the angular dependence is represented: the scale of the graphs is not the same.}
    \label{fig:iclb2d-angle}
\end{figure}

Many additional properties of \(\iclbh\) have been established. Let us introduce \(\calU\defby\setof{(\beta,h)}{h\neq 0, \text{ or } h=0 \text{ and } \beta<\betac(d)}\). We also write, here and in the sequel, \(\uvec_x \defby x/\normII{x}\) for all \(x\in\Rd\setminus\{0\}\).
\begin{itemize}
    \item \(\lim_{\beta\to 0} \iclbh = \lim_{\beta\to \infty} \iclbh = \lim_{h\to \infty} \iclbh = \infty\) and \(\lim_{\beta\uparrow\betac(d)} \iclb = 0\).
    \item For all \((\beta,h)\in\calU\), \(\iclbh(\uvec)\) depends analytically on the direction \(\uvec\)~\cite{Campanino+Ioffe+Velenik-2003,Ott2018}.
    \item For all \((\beta,h)\in\calU\), the function \(\iclbh\) can be extended to a norm on \(\Rd\) by setting \(\iclbh(x)\defby\iclbh(\uvec_x)\normII{x}\) for all \(x\in\Rd\setminus\{0\}\) and \(\iclbh(0)\defby 0\).
    \item For all \((\beta,h)\in\calU\), \(\iclbh\) satisfies the following \define{sharp triangle inequality}~\cite{Pfister+Velenik-1999,Campanino+Ioffe+Velenik-2003,Ott2018}: there exists \(\kappa>0\) such that, for all \(x,y\in\Rd\),
    \[
    \iclbh(x) + \iclbh(y) - \iclbh(x+y) \geq \kappa \bigl( \normII{x} + \normII{y} - \normII{x+y} \bigr).
    \]
    \item When \(d=2\) and \(h=0\), \(\iclbh\) has been computed explicitly (see~\cite{McCoy+Wu-2014} and Fig.~\ref{fig:iclb2d-beta} and~\ref{fig:iclb2d-angle}). 
\end{itemize}


\bigskip

It follows from the definition of the inverse correlation length that, for all \(x\in\Zd\),
\[
\mu_{\beta,h}^+(\sigma_0;\sigma_x) = \exp\bigl[-\bigl(\iclbh(\uvec_x)+\smo(1) \bigr) \normII{x}\bigr],
\]
as \(\normII{x}\to\infty\).
 
Below, the main questions we want to address are the following: what can be said about the corrections to this exponential decay, that is, what is hidden in the \(\smo(1)\) term above? What about more general pairs of local functions?

\section{Coupling to an effective random walk}
\label{sec:OZcoupling}
The nonperturbative approach to sharp asymptotics of correlation functions developed during the last two decades is based on a combination of some of the basic ideas introduced by Ornstein and Zernike, various graphical representations, suitable coarse-graining arguments and ideas originating in perfect simulation algorithms. The outcome is a powerful coupling between the combinatorial objects encoding the correlations, which have complex geometrical and statistical properties, and a directed random walk on \(\Zd\), from which one can then extract the desired asymptotics using standard limit theorems of probability theory.

\subsection{Graphical representations}
This short review is not the proper venue to provide a detailed discussion of the many graphical representations of correlations in the Ising model, and we refer to~\cite{Pfister+Velenik-1999,Duminil-Copin-2017,Friedli+Velenik-2017} for more information. Instead, we will content ourselves with some general comments on the structure such representations take and on some of their basic properties. This will be sufficient to explain the main ideas of the (modern) OZ theory, as well as for the heuristic derivations in later sections.

There are many graphical representations of the Ising model, among which the high- and low-temperature representations, the Fortuin--Kastelyn (also known as FK percolation or random-cluster) representation and the random-current representation. Each has its own advantages and drawbacks. Combined, they yield powerful tools for the rigorous study of the Ising model.

For simplicity, we only present an (informal) discussion of the high-temperature representation; the other ones provide similar expressions, with paths replaced by other, more complicated, connected subsets of the lattice (this higher complexity being compensated, however, by an increased flexibility).

Let us consider the expectation \(\mu_{\beta,0}(\sigma_A)\) where \(A\subset\Zd\) contains \(2n\) vertices and \(\beta<\betac(d)\). The high-temperature representation states that
\begin{equation}\label{GR:HT}
\mu_{\beta,0}(\sigma_A)
=
\sum_{\text{pairings }\pi} \sum_{\substack{\gamma_1,\dots,\gamma_n\sim\pi}} \sfw_\beta(\gamma_1,\dots,\gamma_n),
\end{equation}
where the first sum is over all possible pairings \(\pi\) of the vertices in \(A\) and the second one is over collections of edge-disjoint edge-selfavoiding paths in \(\Zd\), the path \(\gamma_k\) connecting the \(k\)-th pair of vertices in \(\pi\). The nonnegative weight \(\sfw_\beta\) is explicit, but its particular form is of no importance for the present discussion (see~\cite{Pfister+Velenik-1999}, where many useful properties of \(\sfw_\beta\) are established). It suffices to say that it includes complicated infinite-range interactions between the paths (as well as infinite-range self-interactions). Fortunately, these interactions decay fast with the distance, making it possible to establish crucial exponential mixing properties that enable the construction of a coupling with an effective random walk, as described in the next section.

\subsection{Coupling with an effective random walk}
Let us now explain what the modern OZ theory tells us about this type of graphical expansions. Here, we'll further restrict our attention the 2-point function \(\mu_{\beta,0}(\sigma_0\sigma_x)\). In this case, \eqref{GR:HT} becomes
\begin{equation}\label{eq:GeomRepr2PTF}
\mu_{\beta,0}(\sigma_0\sigma_x)
=
\sum_{\gamma:0\to x} \sfw_\beta(\gamma),
\end{equation}
where the sum runs over edge-selfavoiding paths connecting \(0\) and \(x\). The general idea is that the path \(\gamma\) can be decomposed as a concatenation of suitably defined sub-paths, with the intuitive picture that the path \(\gamma\) vehicles the ``full influence'' of the spin at \(0\) on the spin at \(x\), while the sub-paths are the vessels of the ``direct influence'' in the OZ heuristics. To achieve this, one needs to construct sub-paths and the associated weights so that \(\sfw_\beta(\gamma)\) factors along the sub-paths.

\begin{figure}[h]
    \centering
    \resizebox{10cm}{!}{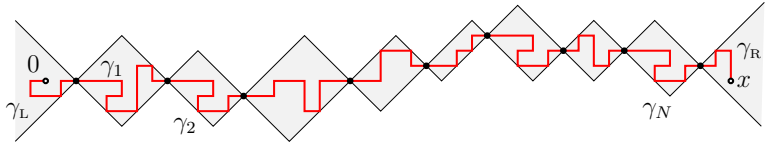}
    \caption{The decomposition of the path \(\gamma\) as a concatenation of \(N\) diamond-confined sub-paths (plus two cone-confined boundary pieces).}
    \label{fig:cone-confinement}
\end{figure}

The first stage of the OZ theory, based on a suitable coarse-graining of the path \(\gamma\) and decoupling ideas first introduced in the context of perfect simulation algorithms, makes this picture rigorous for the following choice of sub-paths: the diamond-confined ones. More precisely, let us introduce the \define{cones} \(\fcone \defby \setof{y\in\Rd}{y\cdot \uvec_x \geq \norm{y}/\sqrt{2}}\), \(\bcone \defby - \fcone\) and the \define{diamonds} \(D(v)\defby \fcone \cap (v+\bcone)\). The sub-paths are then the paths contained in translates of \(D\) and the output of the construction is the identity
\[
\bigl( 1+\sfO(e^{-c \norm{x}}) \bigr)\, e^{\iclb(x)}\sum_{\gamma:0\to x} \sfw_\beta(\gamma)
=
\sum_{\substack{\gamma:0\to x \\ \gamma = \gamma_{\scriptscriptstyle\rm L} \sqcup \gamma_1 \sqcup \cdots \sqcup \gamma_N \sqcup \gamma_{\scriptscriptstyle\rm R}}}
\sfq_{\scriptscriptstyle\rm L}(\gamma_{\scriptscriptstyle\rm L}) \sfq_{\scriptscriptstyle\rm R}(\gamma_{\scriptscriptstyle\rm R}) \prod_{k=1}^N \sfq(\gamma_k) ,
\]
where \(c\) is some positive constant (independent of \(x\)), \(\sfq_{\scriptscriptstyle\rm L}\), \(\sfq_{\scriptscriptstyle\rm R}\) and \(\sfq\) are suitable non-negative weights (depending on \(\uvec_x\), but not on \(\norm{x}\)) and the sum in the right-hand side is taken over paths \(\gamma\) that are concatenation of sub-paths \(\gamma_{\scriptscriptstyle\rm L}\), \(\gamma_1,\dots,\gamma_N\) and \(\gamma_{\scriptscriptstyle\rm R}\) (where \(N\geq 0\) is not fixed) satisfying the following confinement properties: \(\gamma_{\scriptscriptstyle\rm L}\), \(\gamma_1,\dots,\gamma_N\) and \(\gamma_{\scriptscriptstyle\rm R}\) are respectively contained in translates of \(\bcone\), \(D\) and \(\fcone\), as depicted in Figure~\ref{fig:cone-confinement}. The weights \(\sfq_{\scriptscriptstyle\rm L}\), \(\sfq_{\scriptscriptstyle\rm R}\) and \(\sfq\) all decay exponentially with the diameter of the sub-path. As a consequence, all the sub-paths are microscopic and the number \(N\) above can be assumed to be proportional to \(\norm{x}\). 

The introduction of the factor \(e^{\iclb(x)}\) makes the expression in the left-hand side sub-exponential in \(\norm{x}\). In the right-hand side, this factor is implicitly contained in the weights and has, in particular, a crucial impact on the weight \(\sfq\) (see below). 

The factor \(\bigl( 1+\sfO(e^{-c \norm{x}}) \bigr)\) in the right-hand side has two sources: the first one is the existence of atypical paths that cannot be split into sub-paths; the second source is the factorization of weights (all sub-paths of a path \(\gamma\) interact and obtaining a product form as above comes at a cost).

\begin{figure}[h]
    \centering
    \resizebox{10cm}{!}{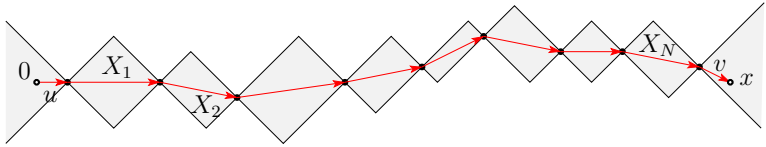}
    \caption{The effective random walk coupled to the microscopic path \(\gamma\).}
    \label{fig:cone-confinement-2}
\end{figure}

It is then convenient to forget about the microscopic path \(\gamma\) and only consider the induced measure on \(\Zd\) defined by the weights (for \(y\in\Zd\))
\[
\rho_{\scriptscriptstyle\rm L}(y) \defby \sum_{\substack{\gamma:0\to y\\\gamma\subset y+\bcone}} \sfq(\gamma),
\qquad
\rho(y) \defby \sum_{\substack{\gamma:0\to y\\\gamma\subset D(y)}} \sfq(\gamma),
\qquad
\rho_{\scriptscriptstyle\rm R}(y) \defby \sum_{\substack{\gamma:0\to y\\\gamma\subset \fcone}} \sfq(\gamma).
\]
All these measures inherit exponential tails.
Moreover, it turns out that the factor \(e^{\iclb(x)}\) that we introduced above automatically guarantees that \(\rho\) is a probability measure. The latter can thus be interpreted as the transition probabilities of a (directed) random walk \(S=(S_k)_{k\geq 0}\) on \(\Zd\).
In terms of these quantities, the outcome of the OZ theory is the following identity:
\begin{multline}\label{eq:OZ2PTF}
\bigl( 1+\sfO(e^{-c \norm{x}}) \bigr)\, e^{\iclb(x)} \sum_{\gamma:0\to x} \sfw_\beta(\gamma)\\
=
\sum_{u,v\in\Zd} \rho_{\scriptscriptstyle\rm L}(u) \rho_{\scriptscriptstyle\rm R}(v) \RWP_u(\exists N\geq 1:\, S_N = x-v),
\end{multline}
where \(\RWP_u\) is the law of the random walk \(S=(S_k)_{k\geq 0}\) starting at \(u\in\Zd\) (see Figure~\ref{fig:cone-confinement-2}). The increments \(X_k\defby S_k-S_{k-1}\) enjoy the following properties:
\begin{itemize}
    \item directedness: \(X_k\cdot \uvec_x > 0\) for all \(k\geq 1\);
    \item exponential tails: \(\exists c>0\) such that \(\RWE_u(e^{c\norm{X_k}}) < \infty\).
\end{itemize}
Moreover, the boundary weights \(\rho_{\scriptscriptstyle\rm L},\,\rho_{\scriptscriptstyle\rm R}\) also exhibit exponential decay: \(\exists c>0\) such that
\[
\rho_{\scriptscriptstyle\rm L}(u)\leq e^{-c\norm{u}},\qquad\rho_{\scriptscriptstyle\rm R}(v) \leq e^{-c\norm{u}}.
\]

It may be interesting at this stage to stress how close the resulting structure is to the OZ heuristics described at the beginning of the introduction.
To do that, let us introduce
\begin{align*}
G(y)
&\defby
\sum_{u,v\in\Zd} \rho_{\scriptscriptstyle\rm L}(u) \rho_{\scriptscriptstyle\rm R}(v) \RWP_u(\exists N\geq 1:\, S_N = x-v)\\
&=
\sum_{u,v\in\Zd} \rho_{\scriptscriptstyle\rm L}(u) \rho_{\scriptscriptstyle\rm R}(v) \sum_{\substack{N\geq 0, y_1,\dots,y_N\\y_1+\cdots+y_N=x-u-v}}\prod_{k=1}^N \rho(y)\\
\end{align*}
and the generating functions (for \(z\in\bbC\))
\begin{alignat*}{4}
\bbG(z) &\defby& \sum_{y\in\Zd} G(y) z^{y\cdot \uvec_x},
\qquad
\bbC(z) &\defby& \sum_{y\in\Zd} \rho(y) z^{y\cdot \uvec_x},\\
\bbB_{\scriptscriptstyle\rm L}(z) &\defby& \sum_{y\in\Zd} \rho_{\scriptscriptstyle\rm L}(y) z^{y\cdot \uvec_x},
\qquad
\bbB_{\scriptscriptstyle\rm R}(z) &\defby& \sum_{y\in\Zd} \rho_{\scriptscriptstyle\rm R}(y) z^{y\cdot \uvec_x}.
\end{alignat*}
A straightforward computation then yields the relation
\[
\bbG(z) = \frac{\bbB_{\scriptscriptstyle\rm L}(z)\bbB_{\scriptscriptstyle\rm R}(z)}{1-\bbC(z)},
\]
which is the analogue of~\eqref{eq:OZEqn}. The appearance of a nontrivial numerator is a consequence of the presence of the ``boundary pieces''. When deriving the asymptotic behavior of \(\mu_{\beta,0}(\sigma_0\sigma_x)\) as \(\norm{x}\) diverges, this will turn out only to change the prefactor by a multiplicative constant.

\subsection{Historical remarks}

In this section, we briefly highlight some of the milestones in the development of the probabilistic (nonperturbative) OZ theory. Note that the development of perturbative approaches (which apply, in principle, to a larger class of systems) is still ongoing, but will not be discussed here.

In the 1980s, inspired by the ideas developed by Ornstein and Zernike, mathematical physicists started to develop a version of the theory allowing a non-perturbative analysis of the asymptotic behavior of correlations under an assumption of finite correlation length, for models in which these correlations can be expressed in terms of sums over suitable classes of connected objects: paths, clusters, etc. The first works dealt with the self-avoiding walk~\cite{Chayes+Chayes-1986}, models of random surfaces~\cite{Abraham+Chayes+Chayes-1985,Jonsson-1986} and Bernoulli percolation~\cite{Campanino+Chayes+Chayes-1991}. Due the approach used, these works suffered from severe limitations: it was crucial that the weights of the geometrical objects enjoy nice factorization properties; the analysis was mostly limited to the axes directions; the proof of separation of masses  was intricate and lacked flexibility.

In the late 1990s and early 2000s, the modern version of the theory was developed in several stages. First, the results for the selfavoiding walk were extended to general directions by importing measure-tilting arguments from the theory of large deviations~\cite{Ioffe-1998}. Then, an easier, much more robust and more powerful way of establishing separation of masses, relying on a suitable coarse-graining of the microscopic object, was introduced in~\cite{Campanino+Ioffe-2002} and applied to Bernoulli percolation (in general directions). The next step was to remove the requirement of factorizability of the weights, thus opening up the possibility of analyzing genuine lattice spin systems. This was achieved in~\cite{Campanino+Ioffe+Velenik-2003} by exploiting exponential mixing in order to couple the original object with a process with infinite memory (described by a suitable Ruelle transfer operator). The first application was to the Ising model above the critical temperature, but the approach was later extended to general Potts model in the same regime~\cite{Campanino+Ioffe+Velenik-2008}. The fact that the approximating process enjoyed suitable mixing properties made it possible to establish local limit theorems for the latter that were generalizations of those available for directed random walks. However, this process remained substantially more difficult to analyze than a random walk, which prevented the applications of these methods to more complicated settings. The last stage was achieved in~\cite{Ott+Velenik-2017}, where it was shown how the original graphical object can be coupled to a genuine (directed) random walk; this was done by adapting decoupling techniques first introduced in the context of perfect simulation in~\cite{Comets+Fernandez+Ferrari-2002}.

\section{Power-law corrections}
\label{sec:powerLawCorrections}
In this section, we review what is known about the corrections to the exponential decay. It turns out that more general covariances have been treated in some cases, so we first discuss this issue.

\subsection{General covariances}

Let \(f\) and \(g\) be two local functions. It is well-known (see~\cite{Friedli+Velenik-2017} for instance) that such functions can be expressed as (finite) linear combinations of random variables of the form \(\sigma_A \defby \prod_{i\in A}\sigma_i\) with \(A\Subset\Zd\). Namely, there exist coefficients \((\hat f_A)_{A\subset\supp(f)}\) and \((\hat g_B)_{B\subset\supp(g)}\) such that
\[
f = \sum_{A\subset\supp(f)} \hat f_A \sigma_A,
\qquad\text{ and }\qquad
g = \sum_{B\subset\supp(g)} \hat g_B \sigma_B.
\]
Using these expressions, we can rewrite
\[
\mu(f;g) = \sum_{\substack{A\subset\supp(f)\\B\subset\supp(g)}} \hat f_A  \hat g_B \, \mu(\sigma_A;\sigma_B) .
\]
This shows that the analysis is essentially reduced to analyzing the asymptotic behavior of the covariances \(\mu (\sigma_A;\sigma_B)\).
\begin{remark}
Note that the derivation of the asymptotic behavior of \(\mu(f;g)\) from the knowledge of the asymptotics of \(\mu (\sigma_A;\sigma_B)\) can only be done for generic functions \(f\) and \(g\)). Indeed, it is always possible to construct pairs of functions (that is, find suitable coefficients \(\hat{f}_A\) and \(\hat{g}_B\)) for which the leading term in the prefactor vanishes.
\end{remark}

\subsection{The regime \(h=0, \beta\in \bigl(0,\betac(d)\bigr)\)}

This is by far the best understood situation. As mentioned above, there is a unique Gibbs measure \(\mu_{\beta,0}\) in this regime. Note that the measure \(\mu_{\beta,0}\) is invariant under a global spin flip, that is, under the mapping \(\omega\mapsto-\omega\) interchanging \(+\) and \(-\) spins. This implies in particular that
\(\mu_{\beta,0}(\sigma_C) = 0\) for all sets \(C\Subset\Zd\) with \(\abs{C}\) odd. This immediately entails that the covariances \(\mu_{\beta,0}(\sigma_A;\sigma_B)\) vanish whenever \(\abs{A}+\abs{B}\) is odd, so that we are left with only two classes of nontrivial covariances: the \define{odd--odd correlations} in which \(\abs{A}\) and \(\abs{B}\) are both odd, and the \define{even--even correlations} in which \(\abs{A}\) and \(\abs{B}\) are both even. It turns out that these two classes display very different asymptotic behavior.

\subsubsection{Odd--odd correlations}
Combining~\eqref{eq:GeomRepr2PTF} and~\eqref{eq:OZ2PTF}, it is not difficult to derive Ornstein--Zernike 
asymptotics for the 2-point function. This can however be extended to general odd-odd correlations.
The best nonperturbative result to date is due to Campanino, Ioffe and Velenik~\cite{Campanino+Ioffe+Velenik-2004} and states that  all odd--odd correlations exhibit OZ decay: for any \(A,B\Subset\Zd\) with \(\abs{A}\) and \(\abs{B}\) odd, there exists an analytic function \(\uvec\mapsto\Psi_{A,B,\beta}(\uvec)\) such that, uniformly in the unit vector \(\uvec\),
\[
\mu_{\beta,0}(\sigma_A;\sigma_{B+[n\uvec]}) = \bigl( 1 + \smo(1) \bigr) \frac{\Psi_{A,B,\beta}(\uvec)}{n^{(d-1)/2}} \, e^{-\iclb(\uvec) n} ,
\]
as \(n\to\infty\).

\medskip\noindent
\textsl{Historical remarks.}
This result has, of course, a long history, going back more than a century ago~\cite{Zernike-1916}. As far as rigorous results are concerned, the 2-point function was studied first. The earliest results relied on exact computations for the planar model~\cite{Wu-1966,McCoy+Wu-2014}. The first results in arbitrary dimensions, but restricted to a perturbative regime (that is, assuming \(\beta \ll 1\)), were proved in~\cite{Abraham+Kunz-1977,Paes-Leme-1978}. Extensions to general odd--odd correlations (still for \(\beta\ll 1\)) were obtained in~\cite{Bricmont+Frohlich-1985a,Bricmont+Frohlich-1985b,Zhizhina+Minlos-1988,Minlos+Zhizhina-1996}. The first nonperturbative proof for the 2-point function was obtained in~\cite{Campanino+Ioffe+Velenik-2003} (and was later extended to Potts models in~\cite{Campanino+Ioffe+Velenik-2008}).

\begin{figure}[h]
    \centering
    \resizebox{6cm}{!}{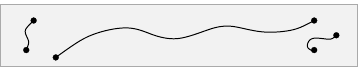}%
    \hspace*{3mm}
    \resizebox{6cm}{!}{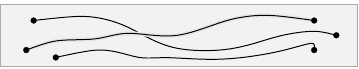}%
    \caption{The two classes of paths configurations: the contribution coming from those with more than one crossing path (right) is negligible.}
    \label{fig:2ptf-HT}
\end{figure}

\medskip\noindent
\textsl{Heuristic derivation.}
Let us now explain heuristically how this result can be obtained. For ease of notation, we only discuss the case \(A=\{x_1,x_2,x_3\}\), \(B=\{y_1,y_2,y_3\}\). The asymptotic behavior of \(\mu_{\beta,0}(\sigma_A;\sigma_{B+[n\uvec]})\) can be understood in various ways. We are going to present a simple one, based on the high-temperature representation of the correlation functions. As discussed above, the latter allows one to write
\[
\mu_{\beta,0}(\sigma_A;\sigma_{B+[n\uvec]}) =
\sum_{\gamma_1,\gamma_2,\gamma_3} \sfw_\beta(\gamma_1,\gamma_2,\gamma_3),
\]
where the sum is over edge-disjoint selfavoiding paths \(\gamma_1\), \(\gamma_2\) and \(\gamma_3\) connecting distinct pairs of vertices in \(\{x_1,x_2,x_3,y_1,y_2,y_3\}\). There are two classes of contributions to the above sum: in the first class, only one of the paths connects some vertex in \(A\) to some vertex in \(B\); in the second one, all three paths do (see Fig.~\ref{fig:2ptf-HT}).

It can be shown that the effective interaction between these paths is \emph{repulsive} (see, for instance, \cite{Pfister+Velenik-1999}). This allows one to get an upper bound of order \(e^{-(3+\sfo(1))\iclb(\uvec)n}\) on the second class of contributions by factorizing the weights of the three paths. This is clearly negligible in view of our target estimate. We can therefore restrict our attention to the first class. Using the coupling described in Section~\ref{sec:OZcoupling}, the determination of the prefactor can be reduced to estimating the probability
\[
\RWP_0(\exists N\,:\, S_N = [n\uvec]) 
\]
that the associated effective directed random walk \(X_k\), starting at \(0\), steps on the vertex \([n\uvec]\).
(Of course, in a rigorous discussion, one should take into account the presence of the other two short paths, as well as the fact that the original path goes from a vertex of \(A\) to a vertex of \(B\), rather than from \(0\) to \([n\uvec]\); however, these issues can be shown to affect only the value of the constant in the prefactor and we thus ignore them here.)

The desired asymptotics now follows from the local CLT for this directed random walk, as we explain now. Let \(\bfD\) be the covariance matrix of the random walk and write \(\bbE(X)=\xi\uvec\) for the expectation of its increments. Fix \(\epsilon>0\). By the local CLT~\cite{Ioffe-1998}, uniformly in \(y\in\Zd\) such that \(\norm{y-k\xi\uvec} < k^{2/3 - \epsilon}\),
\[
\RWP_0(S_k = y) = \frac{1+\sfo_k(1)}{(s\pi k)^{d/2}\sqrt{\det \bfD}} \exp\Bigl( - \frac{\bigl( y-k\xi\uvec, \bfD(y-k\xi\uvec) \bigr)}{2k} \Bigr) ,
\]
where \((x,y)\) denotes the scalar product of \(x,y\in\Rd\).
In our application, the number of steps \(k\) is allowed to vary. However, simple large deviations estimates imply the existence of \(c>0\) such that, uniformly in \(y\) satisfying \(\norm{y-k\xi\uvec} < k^{1/2-\epsilon}\),
\[
\RWP_0(\exists N\,:\, S_N = y, \abs{N-k} > k^{1/2+\epsilon} ) \leq e^{-ck^{2\epsilon}}.
\]
In our case, \(y=[n\xi\uvec]\). From the above, we thus obtain that, up to corrections decaying faster than any power law,
\[
\RWP_0(\exists N\,:\, S_N = [n\uvec])
=
\sum_{N=N_-}^{N_+} \RWP_0(S_N = [n\uvec]),
\]
where we have introduced \(N_\pm \defby (n/\xi) \pm n^{1/2+\epsilon}\). A little algebra then yields
\begin{align*}
\RWP_0(\exists N\,:\, S_N = [n\uvec])
&=
\frac{1+\sfo_n(1)}{(2\pi n/\xi)^{d/2} \sqrt{\det\bfD}} \sum_{\ell = -n^{1/2+\epsilon}}^{n^{1/2+\epsilon}} 
\exp\Bigl( - \frac{\ell^2 \xi^3 (\uvec, \bfD \uvec)}{2n} \Bigr) \\
&=
\frac{1+\sfo_n(1)} {(2\pi n/\xi)^{d/2} \sqrt{\det\bfD}}\, n^{1/2} \int_{-\infty}^{\infty} 
\exp\Bigl( - \frac{x^2 \xi^3 (\uvec, \bfD \uvec)}{2} \Bigr) \, \dd x\\
&=
\frac{1+\sfo_n(1)} {(2\pi/\xi)^{(d-1)/2} \xi \sqrt{\det\bfD} \sqrt{(\uvec,\bfD\uvec)}}\, n^{-(d-1)/2} ,
\end{align*}
which is the desired OZ prefactor. This computation shows that the OZ decay can be interpreted as a CLT-type result in disguise, thereby elucidating the origin of its universal character.

\subsubsection{Even--even correlations}
\label{sec:EvenEvenHeuristic}
The best nonperturbative result to date is due to Ott and Velenik~\cite{Ott+Velenik-2018} and states the following: for any \(A,B\Subset\Zd\) with \(\abs{A}\), \(\abs{B}\) even and any unit vector \(\uvec\), there exist constants \(0<C_-<C_+<\infty\) such that
\[
\frac{C_-}{\Xi_d(n)} \,  \, e^{-2\iclb(\uvec) n}
\leq
\mu_{\beta,0}(\sigma_A;\sigma_{B+[n\uvec]})
\leq
\frac{C_+}{\Xi_d(n)} \,  \, e^{-2\iclb(\uvec) n},
\]
where
\[
\Xi_d(n) \defby
\begin{cases}
n^2				&	\text{when } d=2,\\
(n\log n)^2		&	\text{when } d=3,\\
n^{d-1}			&	\text{when } d\geq 4.
\end{cases}
\]
Note that the rate of exponential decay is twice the one for odd--odd correlations and that the OZ decay never holds for even--even correlations.

\medskip\noindent
\textsl{Historical remarks.}
This result also has a long history. Early results dealt with the energy-energy correlations, that is, the case where both \(A\) and \(B\) consist of two neighboring vertices. First, exact computations for the planar model established the behavior of the prefactor when \(d=2\)~\cite{Stephenson-1966,Hecht-1967}. In the (nonrigorous) physics literature, Polyakov~\cite{Polyakov-1969} predicted the correct behavior, described above, but a conflicting prediction was put forward by Camp and Fisher~\cite{Camp+Fisher-1971} that the prefactor would be of order \(n^{-d}\) in~ every dimension \(d\geq 2\). This debate was first settled in~\cite{Bricmont+Frohlich-1985a,Bricmont+Frohlich-1985b}, in which the correct prefactor was derived rigorously in dimensions \(d\geq 4\) under the assumption that \(\beta\ll 1\). The corresponding derivation including dimensions \(2\) and \(3\) was first done (still assuming \(\beta\ll 1\)) in~\cite{Zhizhina+Minlos-1988,Minlos+Zhizhina-1996}; the latter works also considered general even--even correlations when \(\beta\ll 1\).

\medskip\noindent
\textsl{Heuristic derivation.}
For simplicity, we consider the case \(A=\{x,y\}\) and \(B+n\uvec = \{u,v\}\). Applying the high-temperature representation to the (non truncated) 4-point function \(\mu_{\beta,0}(\sigma_x\sigma_y\sigma_u\sigma_v)\) yields three possible pairings of the four vertices:

\begin{figure}[h]
    \centering
    \resizebox{\textwidth}{!}{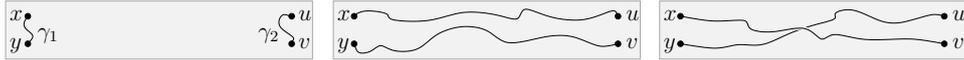}%
    \caption{The three possible pairings of \(x,y,u,v\). }
    \label{fig:eveneven}
\end{figure}%
\noindent
Let us look at the left-most picture. Since \(\beta<\betac(d)\), the paths tend to stay close to their endpoints. In particular, one might expect that, when \(n\gg 1\), \(\sfw_\beta(\gamma_1,\gamma_2) \approx \sfw_\beta(\gamma_1)\sfw_\beta(\gamma_2)\) and that the contribution of such paths should approximately yield
\[
\sum_{\substack{\gamma_1:x\to y\\\gamma_2:u\to v}} \sfw_\beta(\gamma_1,\gamma_2)
\approx
\sum_{\gamma_1:x\to y} \sfw_\beta(\gamma_1) \sum_{\gamma_2:u\to v} \sfw_\beta(\gamma_2)
=
\mu_{\beta,0}(\sigma_x\sigma_y)\mu_{\beta,0}(\sigma_u\sigma_v).
\]
In particular, this would imply that the \emph{truncated} correlation \(\mu_{\beta,0}(\sigma_x\sigma_y;\sigma_u\sigma_v)\) should be approximately given by the contribution coming from pairs of paths crossing from \(A\) to \(B\) as in the two right-most pictures above.

Now, assuming further (rather optimistically) that it is still possible to write \(\sfw_\beta(\gamma_1,\gamma_2) \approx \sfw_\beta(\gamma_1)\sfw_\beta(\gamma_2)\) for such pairs, we would obtain
\begin{align*}
\mu_{\beta,0}(\sigma_x\sigma_y;\sigma_u\sigma_v)
&\approx
\sum_{\gamma_1:x\to u} \sfw_\beta(\gamma_1) \sum_{\gamma_2:y\to v} \sfw_\beta(\gamma_2)
+
\sum_{\gamma_1:x\to v} \sfw_\beta(\gamma_1) \sum_{\gamma_2:y\to u} \sfw_\beta(\gamma_2)\\
&=
\mu_{\beta,0}(\sigma_x\sigma_u)\mu_{\beta,0}(\sigma_y\sigma_v) + \mu_{\beta,0}(\sigma_x\sigma_v)\mu_{\beta,0}(\sigma_y\sigma_u).
\end{align*}
The previously derived OZ asymptotics for 2-point functions would then imply a decay of order \(n^{-(d-1)} e^{-2\iclb(\uvec)}\). This is the correct behavior when \(d\geq 4\), but the approximations made along the way are too coarse to yield the proper order of decay in lower dimensions.

So, let's move back one step and assume that one can factorize \(\sfw_\beta(\gamma_1,\gamma_2) \approx \sfw_\beta(\gamma_1)\sfw_\beta(\gamma_2)\), but that one cannot ignore the constraint that \(\gamma_1\) and \(\gamma_2\) must be edge-disjoint. One would then get an additional correction coming from the probability that two such paths do not intersect. If this probability behaved in the same way as it does for two directed random walk bridges, then this additional factor would be of order
\(n^{-1}\) when \(d=2\), \((\log n)^{-2}\) when \(d = 3\) and \(1\) when \(d\geq 4\). Combining this with the previously obtained prefactor \(n^{-(d-1)}\) then yields the correct behavior.

Let us stress that, although it leads to the correct conclusion, this heuristic argument is not very convincing. The two most glaring
problems are the following. The first one occurs when we approximated the contribution from paths connecting \(x\) to \(y\) and \(u\) to \(v\) by the product of the 2-point functions (thus conveniently canceling the product appearing in the truncation). However, while it is true that \(\sfw_\beta(\gamma_1,\gamma_2) \approx \sfw_\beta(\gamma_1)\sfw_\beta(\gamma_2)\) in that case, the error made is of order \(e^{-(2+\sfo(1))\iclb(\uvec)}\), that is, it is of the same order as our target estimate and thus cannot be neglected. The second problem comes from the factorization \(\sfw_\beta(\gamma_1,\gamma_2) \approx \sfw_\beta(\gamma_1)\sfw_\beta(\gamma_2)\) in configurations with two crossing paths. Indeed, given that the two paths travel together over a long distance, such a claim is by no means obvious. In fact, if the effective interaction between the two paths was \emph{attractive}, then they might ``stick'' to each other and travel together (pinning effect), which would both make the rate of exponential decay strictly smaller and restore standard OZ behavior for the prefactor. In fact, it is precisely this type of problems that still prevents an extension of such an analysis to Potts models (in particular, for these models, we do not even know whether the rate of exponential decay is twice the inverse correlation length or is strictly smaller!).

\begin{figure}
    \centering
    \includegraphics[width=\textwidth]{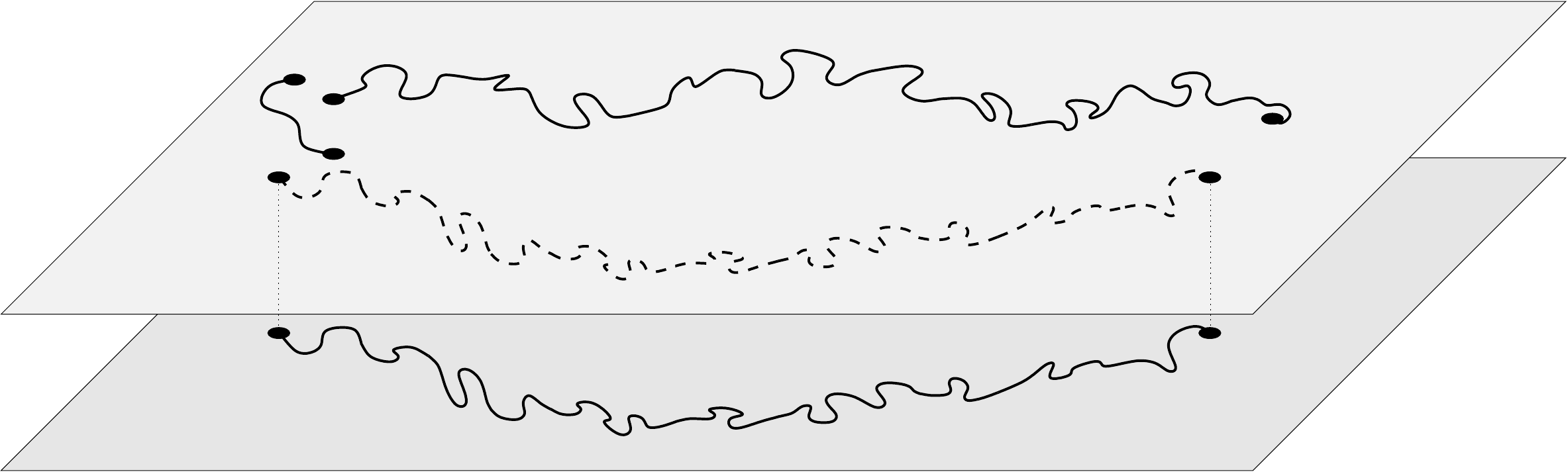}
    \caption{The truncated correlation function \(\mu_{\beta,0}(\sigma_A;\sigma_B)\) with \(\abs{A}=4\) and \(\abs{B}=2\) can be bounded from above in terms of a pair of crossing high-temperature paths living in independent copies of the system and conditioned not to intersect. A similar expression in terms of FK-clusters provides a lower bound.}
    \label{fig:noncrossing}
\end{figure}

In view of this, it is quite remarkable that a version of this argument can actually be made rigorous, at least as upper and lower bounds. First, one obtains an upper bound on \(\mu_{\beta,0}(\sigma_x\sigma_y;\sigma_u\sigma_v)\) in terms of pairs of crossing high-temperature paths living in \emph{independent} copies of the system, but \emph{constrained not to intersect}; see Figure~\ref{fig:noncrossing}. Similarly, the lower bound is expressed in terms of pairs of crossing FK-clusters living in \emph{independent} copies of the system, but \emph{constrained not to intersect}. We emphasize that these are bounds, not asymptotic identities. In order to obtain asymptotically matching upper and lower bounds, one should probably obtain a good control of the interaction of the two paths and not simply bound it away.

All in all, the above arguments show that \(e^{2\iclb(\uvec)}\mu_{\beta,0}(\sigma_x\sigma_y;\sigma_u\sigma_v)\) is indeed comparable to the probability that two \emph{independent} random walks starting from \(\{x,y\}\) are ending at \(\{u,v\}\) without intersecting.

\medskip\noindent
\textsl{Links to critical behavior.}
Let us conclude this section by emphasizing the nontrivial dimension-dependent behavior occurring in this case, which is reminiscent of what is commonly seen at criticality. (Similarly rich behavior is also observed in many other instances; see, for example, \cite{Ott+Velenik-2017}). This has been discussed heuristically in~\cite{Bricmont+Frohlich-1985a}, from which we extracted the following excerpt:
\begin{quote}\small
    There is a striking similarity between these power-law corrections and the power-law decay at the critical point in at least two respects: there is a critical dimension, equal to three here, above which mean field theory is correct and which is characterized by logarithmic corrections to mean field theory. Moreover, the above discussion shows the relevance of intersection properties of random walks, a typical feature of critical point theory. However, here we deal with intersections of two walks at the same time, while it appears that the intersection of the paths is the relevant quantity for critical phenomena. Although one is dealing with a high-temperature situation where the theory is massive and non-critical, the similarity with critical phenomena is as follows: the high-temperature lines joining 0 and x have a weight that is exponentially decreasing with their length, and this produces the mass. However, the power-law corrections, as our analyses have demonstrated, have their origin in the transverse fluctuations of the random line. These are not massive in any sense, and the question of their Gaussian or non-Gaussian nature is equivalent to the Ornstein- Zernike (or mean field) decay or to violations of this decay.
\end{quote}

\subsection{The regime \(h\neq 0, \beta\in (0,+\infty)\)}

Again, there is a unique Gibbs measure \(\mu_{\beta,h}\). The best nonperturbative result to date is due to Ott~\cite{Ott2018} and states that the 2-point function displays OZ decay: there exists an analytic function \(\uvec\mapsto\Psi_{\beta,h}(\uvec)\) such that, uniformly in the unit vector \(\uvec\),
\[
\mu_{\beta,h}(\sigma_0;\sigma_{[n\uvec]}) = \bigl( 1 + \smo(1) \bigr) \frac{\Psi_{\beta,h}(\uvec)}{n^{(d-1)/2}} \, e^{-\iclbh(\uvec) n} ,
\]
as \(n\to\infty\). One expects that similar arguments should prove that this result extends to \emph{all} truncated functions of the form \(\mu_{\beta,h}(\sigma_A;\sigma_B)\) (that is, even-even correlations are not expected to enjoy non-OZ asymptotics once \(h\neq 0\)).

\begin{figure}[h]
    \centering
    \resizebox{!}{9.5mm}{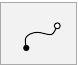}%
    \hspace*{3mm}
    \resizebox{!}{9.5mm}{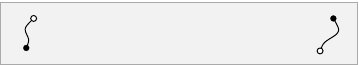}%
    \hspace*{3mm}
    \resizebox{!}{9.5mm}{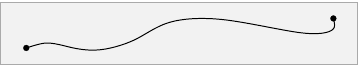}%
    \caption{Left: Paths contributing to \(\mu_{\beta,h}(\sigma_i)\) start at \(i\) and end at an arbitrary vertex (represented in white). Middle and right: The two classes of contributions to \(\mu_{\beta,h}(\sigma_0\sigma_{[n\uvec]})\).}
    \label{fig:2ptfField-HT}
\end{figure}

\medskip\noindent
\textsl{Heuristic derivation.}
As in the case of even--even correlations above, one of the difficulties here is that the truncation in \(\mu_{\beta,h}(\sigma_0;\sigma_{[n\uvec]})\) is non-trivial (that is, \(\mu_{\beta,h}(\sigma_0)\mu_{\beta,h}(\sigma_{[n\uvec]})\neq 0\)). The asymptotic behavior can thus only be seen once the leading contribution from the non-truncated correlation \(\mu_{\beta,h}(\sigma_0\sigma_{[n\uvec]})\) has been canceled. This problem is solved by exploiting properties of the random-current representation. Nevertheless, in this heuristic discussion, we stick to the high-temperature one. In the presence of  magnetic field, the latter states that, for any \(i\in\Zd\),
\[
\mu_{\beta,h}(\sigma_i) = \sum_{\gamma} \sfw_{\beta,h}(\gamma),
\]
where the sum is over finite edge-selfavoiding paths starting at \(i\) and ending anywhere, see Figure~\ref{fig:2ptfField-HT}, left. Similarly,
\[
\mu_{\beta,h}(\sigma_0\sigma_{[n\uvec]}) = \sum_{\gamma_1,\gamma_2} \sfw_{\beta,h}(\gamma_1,\gamma_2) + \sum_{\gamma} \sfw_{\beta,h}(\gamma),
\]
where the first sum is over pairs of edge-disjoint finite edge-selfavoiding paths starting at \(0\) and \([n\uvec]\) respectively, while the second sum is over edge-selfavoiding paths starting at \(0\) and ending at \([n\uvec]\); see Figure~\ref{fig:2ptfField-HT}, middle and right.

As in the case of even--even correlations (and with the same problems when actually implementing this idea), it seems natural to expect that the contributions coming from pairs of paths should factor, as these paths tend to be small. This would again lead to
\[
\mu_{\beta,h}(\sigma_0;\sigma_{[n\uvec]}) \approx \sum_{\gamma} \sfw_{\beta,h}(\gamma),
\]
where the sum is over edge-selfavoiding paths starting at \(0\) and ending at \([n\uvec]\), as in Figure~\ref{fig:2ptfField-HT}, right. We are thus back to a picture resembling the one appearing when analyzing the 2-point function at \(h=0\) and \(\beta<\betac(d)\) (but with a different weight for the path), which makes it plausible that OZ decay should also hold in this case.

The rigorous version of this argument in~\cite{Ott2018} uses the random-current representation and is based on the identity~\eqref{eq:GeomReprTrunc2PTF} derived in Appendix~\ref{app:RCmodel}.

\subsection{The regime \(h=0, \beta\in (\betac(d),+\infty)\)}
\label{sec:LTheuristic}
In this regime, there are multiple extremal Gibbs measures and results are restricted to \(\mu_{\beta,0}^+\) (and \(\mu_{\beta,0}^-\) by symmetry).
This is the regime in which the understanding is most incomplete, even at the level of the truncated 2-point function. 
The only nonperturbative results rely on exact computations and are thus restricted to the planar model. They have been obtained by Wu~\cite{Wu-1966} and state that OZ asymptotics are violated for the truncated 2-point function: there exists an (explicit) constant \(C\in(0,\infty)\) such that
\[
\mu_{\beta,0}(\sigma_0;\sigma_{[n\eone]}) = \bigl( 1 + \smo(1) \bigr) \frac{C}{n^2} \, e^{-\iclb(\eone) n} ,
\]
as \(n\to\infty\). (See~\cite{WuEtAl-1976} for general directions \(\uvec\).)

When \(d\geq 3\), only perturbative results are available. They have been obtained by Bricmont and Fröhlich~\cite{Bricmont+Frohlich-1985b} and state that OZ asymptotics hold for the truncated 2-point function at sufficiently low temperatures: there exist constants \(0<C_-<C_+<\infty\) such that
\[
\frac{C_-}{n^{(d-1)/2}} \,  \, e^{-\iclb(\uvec) n}
\leq
\mu_{\beta,0}(\sigma_0;\sigma_{[n\eone]})
\leq
\frac{C_+}{n^{(d-1)/2}} \,  \, e^{-\iclb(\uvec) n}.
\]

\medskip\noindent
\textsl{Heuristic derivation.} In this regime, it is easier to understand the behavior in terms of the low-temperature representation, that is, in terms of Peierls contours (that is, the connected components of dual plaquettes separating \(+\) and \(-\) spins; see~\cite{Friedli+Velenik-2017}). Namely, one might argue that, roughly speaking, correlations between two distant vertices are due to the presence of a large contour surrounding both vertices. Strictly speaking, this is incorrect: although it can be used to obtain an upper bound on the truncated 2-point function, this bound is not sharp. Nevertheless, it is possible to make sense of such a picture, either in a perturbative regime (\(\beta\gg 1\)), where cluster expansion can be used to express the truncated 2-point function in terms of a ``chain of contours'' surrounding both vertices, or non perturbatively using the random-current representation, as is explained in Appendix~\ref{app:RCmodel}. For simplicity, we pursue this discussion with Peierls contours, as this does not affect the main issues.

\begin{figure}[h]
    \centering
    \resizebox{!}{1.5cm}{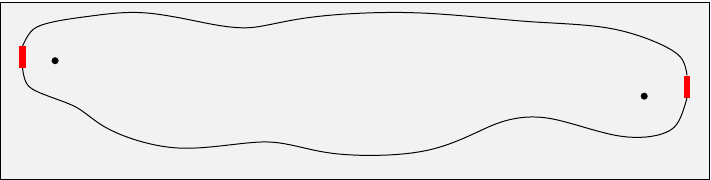}%
    \hspace*{5mm}
    \resizebox{!}{1.5cm}{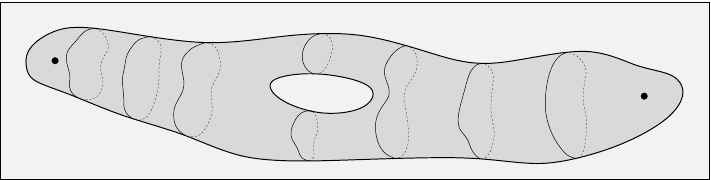}%
    \caption{Left: The planar case; the contour surrounding \(0\) and \([n\uvec]\) can be seen as a pair of edge-disjoint open paths. Right: The 3-dimensional case.}
    \label{fig:2ptf-LT}
\end{figure}

Let us first consider the planar model. In this case (represented in Figure~\ref{fig:2ptf-LT}, left), the contour surrounding the two vertices can be reinterpreted as a pair of edge-disjoint paths. This should be reminiscent of the picture obtained for energy-energy correlations. This makes it plausible that the same type of asymptotics should apply here, which turns out to indeed be the case. So, the interpretation of the anomalous prefactor in the planar case (an exponent equal to \(2\) rather than \(1/2\)) can be attributed to the entropic repulsion between the two ``halves'' of the Peierls contour, which becomes ``fat'' (its width being of order \(n^{1/2}\)).

The situation is very different both when the model is not planar, that is, when \(d\geq 3\) (see Figure~\ref{fig:2ptf-LT}, right) or when \(d=2\) and the interaction is not restricted to nearest-neighbors. In such cases, since the cost of the contour grows proportionally with its surface area, the cost of having a large section becomes too high and the section remains bounded on average, at least at sufficiently low temperatures, and the contour behaves roughly like a string connecting the two vertices. This leads to standard OZ decay, as proved in~\cite{Bricmont+Frohlich-1985b}. However, it is conceivable, although highly unlikely, that a transitions occurs at a temperature below the critical temperature, at which the section of the contours starts to diverge with \(n\). This so-called \define{breathing} transition does occur in some models of random surfaces, but is not expected to occur here. Proving this remains an important open problem.

\subsection{The regime \(h=0, \beta=\betac(d)\)}
Once more, there is a unique Gibbs measure \(\mu_{\betac(d),0}\)~\cite{Aizenman+Duminil-Copin+Sidoravicius-2015}. In this regime, the results are more satisfactory than for \(\beta>\betac(d)\), but still only partial.

In the case of the planar Ising model, one can again rely on the latter's integrability to compute explicitly the asymptotic behavior of the 2-point function. This was done by Wu~\cite{Wu-1966}, whose result states that there exists an (explicit) constant \(C\in(0,\infty)\) such that
\[
\mu_{\betac(2),0}(\sigma_0;\sigma_i) = \frac{C}{\norm{i}^{1/4}} \bigl( 1 + \sfo(1) \bigr) .
\]
When \(d=3\), only rough polynomial bounds are available: there exist \(c_-,c_+\in(0,\infty)\) such that, for all \(i\in\bbZ^d\),
\[
\frac{c_-}{\norm{i}^2} \leq \mu_{\betac(3),0}(\sigma_0;\sigma_i) \leq \frac{c_+}{\norm{i}} .
\]
At the moment of writing, the best numerical estimate of the exponent in the prefactor is \(-1.036298\)~\cite{KosEtAl-2016}.

The situation is better in higher dimensions. First, for all \(d\geq 5\), there exist \(c_-,c_+\in(0,\infty)\) such that, for all \(i\in\bbZ^d\),
\[
\frac{c_-}{\norm{i}^{d-2}} \leq \mu_{\betac(d),0}(\sigma_0;\sigma_i) \leq \frac{c_+}{\norm{i}^{d-2}} .
\]
Even more is known when either \(d>4\) and the interaction is sufficiently spread-out, or \(d\) is large enough and the interaction is nearest-neighbor. Namely, Sakai~\cite{Sakai-2007, Sakai-2020} proved the existence of a constant \(C\in (0,\infty)\) such that 
\[
\mu_{\betac(d),0}(\sigma_0;\sigma_i) = \frac{C}{\norm{i}^{d-2}} \bigl( 1 + \sfo(1) \bigr) .
\]
\begin{remark}
    There has been recent progress~\cite{Slade+Tomberg-2016} showing that the decay is \(\sfO(\norm{i}^{-2})\) in the \(\varphi^4\) model in dimension \(4\). The latter is expected to have the same behavior as the 4-dimensional Ising model.
\end{remark}

\section{Higher Order Correlations and Mixing Properties}

\subsection{Relation to the pressure}

A lot of information on the Ising model is encoded in the \define{pressure} (or \define{free energy}):
\begin{equation*}
	\psi(\beta,h) = \lim_{N\to\infty} \frac{1}{|\Lambda_N|} \log(\sfZ_{\Lambda_N}),
\end{equation*}
where \(\Lambda_N\) is the square box of sidelength \(N\). In particular, phase transition points are usually defined as the set of points where \(\psi(\beta,h)\) fails to be analytic. Higher-order correlations are closely related to the analyticity of the pressure. Indeed, analyticity in \(h\) is implied by a bound of the form (called a \define{tree bound})
\begin{equation*}
	|U_n(\sigma_{x_1},\sigma_{x_2},\dots,\sigma_{x_n})| \leq C^n e^{-m t(x_1,\dots,x_n)}
\end{equation*}
for some \(C\geq 0, m>0\) not depending on \(n\), where \(t(x_1,\dots,x_n)\) is the length of the Steiner tree generated by \(x_1,\dots,x_n\). It is believed to be satisfied in the whole regime where the pressure is analytic in \(h\) (see~\cite{Seiler-1982}). A stronger property than weak mixing (see next section), called \define{strong mixing} (see~\cite{Dobrushin-Shlosman-1987}), is known to hold at sufficiently high temperature and to imply the tree bound. See Section~\ref{section:what_is_known} for further discussion on weaker statements.

\subsection{Exponential Mixing}

The Ursell functions are a way to measure how well different regions of space decouple. Another way to measure this is through the \define{weak mixing property}: there exist \(c>0, C\geq 0\) such that for any \(\Delta, \Delta'\subset \Zd\) and any two events \(A,B\) with support in \(\Delta, \Delta'\) respectively,
\begin{equation*}
	\big|\mu(A\cap B) - \mu(A)\mu(B)\big| \leq C\sum_{i\in \Delta, j\in\Delta'} e^{-c\norm{i-j}}.
\end{equation*}
This property holds whenever the system is not at a transition point (that is, \((\beta,h)\notin [\betac,\infty)\times \{0\}\)); see the next section for references. Whether this property can be deduced from exponential decay of truncated two-point functions and uniqueness of the infinite-volume measure and whether it is sufficient to obtain optimal bounds on Ursell functions are open problems. (Weak mixing is deduced from finite-volume exponential relaxation of the magnetization well inside a box \emph{uniformly over boundary conditions}, which is \textit{a priori} stronger than decay of the truncated two-point functions. \textit{A posteriori}, these conditions are equivalent for the Ising model.)

\subsection{What is Known}
\label{section:what_is_known}

After the first Ursell function, covariances are by far the most studied Ursell functions. Some information is nevertheless available on higher-order correlations. Using correlation inequalities, one is for example able to deduce that
\begin{equation}
\label{eq:UB_Ursell}
	|U_n(\sigma_{x_1},\sigma_{x_2},\dots,\sigma_{x_n})| \leq C_n e^{-m_n t(x_1,\dots,x_n)},
\end{equation}
where \(m_n>0\) and \(C_n\) are functions of \(n\), whenever one has exponential decay of suitable covariances~\cite{Lebowitz-1972,MartinLof-1973} (see also~\cite[section II.12]{Simon-1993}). However, the dependence of \(m_n\) and \(C_n\) on \(n\) is so bad that, while these bounds can give smoothness of the pressure, they do not imply analyticity. In a different direction, at \(h=0\), using the random-current representation, it has been possible to prove that \((-1)^{n+1} U_{2n}(\sigma_{x_1},\sigma_{x_2},\dots,\sigma_{x_{2n}}) \geq 0\)~\cite{Shlosman-1986}. As an example of the usefulness of such correlation inequalities, the \define{GHS inequality}~\cite{Griffiths+Hurst+Sherman-1970}
\begin{equation*}
	U_3(\sigma_{x_1},\sigma_{x_2},\sigma_{x_3})\leq 0
\end{equation*}
for \(h\geq 0\), easily implies continuity of the magnetization in \(h\) on \((0,\infty)\) for any \(\beta\geq 0\).

Regarding analyticity of the pressure, it is known that the pressure is non-analytic if and only if \(h=0, \beta\geq \betac\). See below for the analyticity statement. The non-analyticity at \(h=0,\beta\geq\betac\) follows from the positivity of the magnetization when \(\beta>\betac\). The nature of the singularity at \(h=0\) was investigated by Isakov~\cite{Isakov-1984} and found to be an essential one when \(\beta\) is sufficiently large. Namely he proved that
\begin{equation*}
	\frac{\partial^k}{(\partial h)^k}\psi(\beta,h)\Big|_{h=0^+} \sim C^k (k!)^{\frac{d}{d-1}}.
\end{equation*}
Extension of this result to any \(\beta>\betac\) remains open.

On the side of mixing properties, weak mixing is known to hold when \(\beta<\betac\) or when \(h\neq 0\) and to imply analyticity of the pressure (see~\cite{Ott-2019} for proof of the two later statements and for further discussion about earlier results). Notice that no non-perturbative result on the tree bound is available.

\section{Open Problems}
We list here some open problems related to what is discussed in the previous sections.

\medskip\noindent
\textsl{Inverse correlation length.}
\begin{enumerate}
    \item Determine whether \(\iclbh\) is analytic in \(\beta\) and/or \(h\) when \((\beta,h)\neq (\betac(d),0)\).
    \item Show that \(\icl_{\betac(d),0} = 0\) for general finite-range interactions in any dimensions. (One way would be to extend the identity \(m^*_{\betac(d)} = 0\) to finite-range interactions).
\end{enumerate}

\medskip\noindent
\textsl{Asymptotics of covariances.}
\begin{enumerate}
    \item Prove sharp asymptotics in the following cases:
    \begin{enumerate*}
        \item Even--even correlations when \(h=0,\beta<\betac(d)\).
        \item In dimension \(2\), determine the order of the prefactor when \(\beta>\betac(2),h=0\) for non-planar interactions. It is expected that one recovers OZ behavior in this case (see~\cite[Section~VII]{Bricmont+Frohlich-1985b}).
        \item Arbitrary \(\mu(\sigma_A;\sigma_B)\) when \(h=0,\beta>\betac(d)\), in particular prove the absence of a breathing transition.
        \item Arbitrary \(\mu(\sigma_A;\sigma_B)\) when \(h\neq 0\).
    \end{enumerate*}
    \item Extend all the results to exponentially decaying interactions of infinite range. Note that the situation is rather subtle in that case: whether the 2-point function enjoys OZ decay actually depends on \emph{corrections} to the exponential decay of the interaction; this will be discussed in the forthcoming~\cite{IOV-2020}.
    \item Remove the assumption that the interactions be purely ferromagnetic. In general, this should lead to an oscillatory prefactor.
\end{enumerate}
    
\medskip\noindent
\textsl{Higher-order correlations.}
\begin{enumerate}
    \item Prove nonperturbatively the tree bound on \(U_n(\sigma_{x_1},\dots,\sigma_{x_n})\).
    \item Determine the exact rate of decay of \(U_n(\sigma_{x_1},\dots,\sigma_{x_n})\) as \(\norm{x_i-x_j}\to\infty\) for all \(i\neq j\).
    \item Then, determine the prefactor.
\end{enumerate}

\medskip\noindent
\textsl{Critical 2-point function.}
\begin{enumerate}
    \item Derive the asymptotic behavior when \(d\geq 5\) for the nearest-neighbor model.
    \item Derive the asymptotic behavior when \(d=3\) or \(4\).
\end{enumerate}

\section{Acknowledgments}
The authors gratefully acknowledge the support of the Swiss National Science 
Foundation through the NCCR SwissMAP.

\appendix

\section{The Random Current Model}
\label{app:RCmodel}

In this Appendix, we present briefly a well-known representation of the Ising model in which one can derive rigorous versions of the basic representations of correlation functions used in the heuristics presented in section~\ref{sec:powerLawCorrections}. We discuss it in finite volume (a finite graph \((\Lambda,E_{\Lambda})\), \(\Lambda\Subset\Zd\)) but an infinite-volume version is also available (see \cite{Aizenman+Duminil-Copin+Sidoravicius-2015}).

\subsection{\(\beta<\betac\), \(h=0\)}

In this regime, the Gibbs measure is unique, so we can choose to work with the free boundary condition (which corresponds to setting \(\eta\equiv 0\)). The random-current representation is obtained by Taylor-expanding \(e^{\beta \sigma_i\sigma_j }\): for \(A\subset\Lambda\),

\begin{align*}
\sum_{\sigma\in\{-1,1\}^{\Lambda}} \sigma_A \prod_{ij\in E_{\Lambda}} e^{\beta \sigma_i\sigma_j } &= \sum_{\sigma\in\{-1,1\}^{\Lambda}} \sigma_A\sum_{\bfn:E_{\Lambda}\to\Zplus} \prod_{ij\in E_{\Lambda}} \frac{\beta^{\bfn_{ij}}}{\bfn_{ij}!} \sigma_{i}\sigma_{j}\\
&= \sum_{\bfn:E_{\Lambda}\to\Zplus} \weight(\bfn) \prod_{i\in\Lambda}\Big(\sum_{\sigma_i\in\{-1,1\}} (\sigma_{i})^{I_i(\bfn) +\mathds{1}_{i\in A}}\Big)\\
&= 2^{|V|}\sum_{\bfn\in \calN_{A}} \weight(\bfn),
\end{align*}
where \(\weight(\bfn)\equiv\weight_{\beta,0}(\bfn) \defby \prod_{e\in E_{\Lambda}}\frac{\beta^{\bfn_e}}{\bfn_e!} \), \(I_{i}(\bfn) \defby \sum_{j: ij\in E_{\Lambda}}\bfn_{ij} \) is the \define{incidence} of \(i\) and
\begin{equation*}
\calN_A \equiv \calN_A(\Lambda) \defby \setof{\bfn:E_{\Lambda}\to\Zplus}{\partial\bfn = A}
\end{equation*}
with \(\partial\bfn \defby \setof{i}{I_{i}(\bfn) = 1\mod 2}\) is the set of \define{sources} of \(\bfn\) (the vertices with odd incidence). Notice that \(\calN_A\) is empty if \(A\) contains an odd number of vertices. We will write
\begin{equation*}
\RCPF(A) \equiv \RCPF_{\Lambda}(A) \defby \sum_{\bfn\in \calN_A} \weight(\bfn).
\end{equation*}
We thus have
\begin{equation*}
\langle \sigma_{A}\rangle_{\Lambda,\beta,0} = \frac{\RCPF(A)}{\RCPF(\varnothing)}.
\end{equation*}

One can interpret a \define{current} \(\bfn\) as a percolation configuration on \(E_{\Lambda}\) by looking at \(\hat{\bfn}_e = \mathds{1}_{\bfn_e>0}\). At this stage, we can already see the structure emerge: by partitioning according to the cluster of \(x,y\) in \(\hat{\bfn}\),
\begin{equation}
\label{eq:RC_HT_TwoPtsFunction}
\langle \sigma_{x}\sigma_y\rangle_{\Lambda} = \frac{\RCPF(x,y)}{\RCPF(\varnothing)} = \sum_{C\ni x,y} \sum_{\bfn\in\calN_{\{x,y\}}} \frac{\weight(\bfn)}{\RCPF(\varnothing)}\mathds{1}_{C_{x,y}(\hat{\bfn}) = C} \defby \sum_{C\ni x,y} \sfw(C),
\end{equation}
with the cluster connecting \(x\) and \(y\) playing the role of the weighted path described previously.

What makes this representation particularly powerful is its ability to deal with duplicated systems. Let us introduce the notation
\begin{equation*}
\RCPF(A)\RCPF(B)\{F\} \defby \sum_{\bfn\in\calN_A}\sum_{\bfm\in\calN_B} \weight(\bfn)\weight(\bfm) F(\bfn+\bfm)
\end{equation*}
for any \(A,B\subset \Lambda\) and function \(F\). We can now state the main feature of the random-current representation: the \define{Switching Lemma}.
\begin{lemma}[Switching Lemma]
	\begin{equation*}
	\RCPF(A)\RCPF(B)\{F\} = \RCPF(A\difsym B)\RCPF(\varnothing)\{ \mathds{1}_{\widehat{\bfn+\bfm}\in\evenPart_B}F \},
	\end{equation*}
	where \(\evenPart_B\) is the event that every cluster of \(\widehat{\bfn+\bfm}\) contains an even number of vertices in \(B\) (possibly \(0\)) and \(\difsym\) denotes symmetric difference.
\end{lemma}

As a first example of how this Lemma can be used, we provide a rigorous version of the graphical representation used in the heuristics presented in Section~\ref{sec:powerLawCorrections} for even--even correlations. We consider the particular case of pair--pair correlations: let \(u,v,x,y\in\Lambda\). Then,
\begin{align*}
\langle \sigma_{x}\sigma_y; \sigma_u\sigma_v\rangle_{\Lambda} &= \frac{1}{\RCPF(\varnothing)^2} \Big(\RCPF(u,v,x,y)\RCPF(\varnothing)-\RCPF(u,v)\RCPF(x,y) \Big)\\
&= \frac{1}{\RCPF(\varnothing)^2} \RCPF(u,v,x,y)\RCPF(\varnothing)\{\mathds{1}_{x\nleftrightarrow y}\}.
\end{align*}
Now the source constraint in the first current implies that, on the event \(\{x\nleftrightarrow y\}\), either \(\{x\leftrightarrow u\}\) or \(\{x\leftrightarrow v\}\). Thus,
\begin{multline*}
\RCPF(u,v,x,y)\RCPF(\varnothing)\{\mathds{1}_{x\nleftrightarrow y}\} = \RCPF(u,v,x,y)\RCPF(\varnothing)\{\mathds{1}_{x\nleftrightarrow y}\mathds{1}_{x\leftrightarrow u}\} +\\+ \RCPF(u,v,x,y)\RCPF(\varnothing)\{\mathds{1}_{x\nleftrightarrow y}\mathds{1}_{x\leftrightarrow v}\}.
\end{multline*}
Using again the Switching Lemma, one obtains
\begin{equation*}
\RCPF(u,v,x,y)\RCPF(\varnothing)\{\mathds{1}_{x\nleftrightarrow y}\mathds{1}_{x\leftrightarrow v}\} = \RCPF(u,y)\RCPF(x,v)\{\mathds{1}_{x\nleftrightarrow y}\}
\end{equation*}
and similarly with \(v\) and \(u\) interchanged. Therefore, we have the identity
\begin{equation*}
\langle \sigma_{x}\sigma_y; \sigma_u\sigma_v\rangle_{\Lambda} = \frac{\RCPF(u,y)\RCPF(x,v)\{\mathds{1}_{x\nleftrightarrow y}\}}{\RCPF(\varnothing)^2} + \frac{\RCPF(v,y)\RCPF(x,u)\{\mathds{1}_{x\nleftrightarrow y}\}}{\RCPF(\varnothing)^2}.
\end{equation*}
We then have
\begin{equation}
\label{eq:RC_HT_EneEneCov}
\frac{\RCPF(u,y)\RCPF(x,v)\{\mathds{1}_{x\nleftrightarrow y}\}}{\RCPF(\varnothing)^2} = \sum_{C_1\ni x,v} \sum_{C_2\ni u,y} \sfw_1(C_1)\sfw_2(C_2)\mathds{1}_{C_1\cap C_2=\varnothing} I(C_1,C_2),
\end{equation}
where \(\sfw_1(C) \defby \frac{\RCPF(x,v)\{\mathds{1}_{C_{x,v}}=C\}}{\RCPF(\varnothing)}\), similarly for \(\sfw_2\), are the same weights as in~\eqref{eq:RC_HT_TwoPtsFunction} and \(I(C_1,C_2)\) is the ratio term making \eqref{eq:RC_HT_EneEneCov} an identity.

\subsection{\(\beta\geq \betac\), \(h=0\)}

In this case, the measure is no longer unique and one wants to study \(\mu^+_{\beta,0}\). We thus add a boundary condition to our setting. This can be done by adding a vertex, denoted \(\partial\), to \(\Lambda\) and an edge between \(\partial\) and \(i\in\Lambda\) if \(i\) shares an edge with a vertex in \(\Zd\setminus\Lambda\). We denote this augmented graph \((\Lambda_{\partial},E_{\Lambda_{\partial}})\). The same expansion as in the previous section gives
\begin{equation*}
\sum_{\sigma\in\{-1,1\}^{\Lambda}} \sigma_A  e^{\beta \sum_{ij\in E_{\Lambda}}\sigma_i\sigma_j +\beta\sum_{i\partial\in E_{\Lambda_{\partial}}} \sigma_i } = 2^{|\Lambda|}
\begin{cases}
\RCPF_{\Lambda_{\partial}}(A) & \text{ if } |A| \text{ even},\\
\RCPF_{\Lambda_{\partial}}(A\cup\{\partial\}) & \text{ if } |A| \text{ odd}.
\end{cases}
\end{equation*}
As a direct consequence,
\begin{equation*}
\langle \sigma_{A}\rangle_{\Lambda,\beta,0}^+ = \frac{1}{\RCPF_{\Lambda_{\partial}}(\varnothing)}\begin{cases}
\RCPF_{\Lambda_{\partial}}(A) & \text{ if } |A| \text{ even},\\
\RCPF_{\Lambda_{\partial}}(A\cup\{\partial\}) & \text{ if } |A| \text{ odd}.
\end{cases}
\end{equation*}
Of particular interest is the truncated two-point function:
\begin{equation*}
\langle \sigma_{x};\sigma_y\rangle_{\Lambda}^+ = \frac{1}{\RCPF_{\Lambda_{\partial}}(\varnothing)^2}\Big( \RCPF_{\Lambda_{\partial}}(x,y)\RCPF_{\Lambda_{\partial}}(\varnothing)-\RCPF_{\Lambda_{\partial}}(x,\partial)\RCPF_{\Lambda_{\partial}}(y,\partial) \Big).
\end{equation*}Now, using the Switching Lemma, one gets
\begin{equation*}
\RCPF_{\Lambda_{\partial}}(x,\partial)\RCPF_{\Lambda_{\partial}}(y,\partial) = \RCPF_{\Lambda_{\partial}}(x,y)\RCPF_{\Lambda_{\partial}}(\varnothing) \{\mathds{1}_{y\leftrightarrow \partial}\},
\end{equation*}and thus
\begin{equation}
\label{eq:RC_LT_TruncTwoPtsFunct}
\langle \sigma_{x};\sigma_y\rangle_{\Lambda}^+ = \frac{\RCPF_{\Lambda_{\partial}}(x,y)\RCPF_{\Lambda_{\partial}}(\varnothing)\{\mathds{1}_{y\nleftrightarrow \partial}\} }{\RCPF_{\Lambda_{\partial}}(\varnothing)^2}
=\sum_{C\subset \Lambda_{\partial}} \mathds{1}_{\partial\notin C} \sfw(C).
\end{equation}
The cluster of \(x,y\) plays the role of the weighted path described in Section~\ref{sec:LTheuristic} and the cost for long cluster is exactly the cost to not connected the cluster of \(x,y\) to \(\partial\) via the background. The weight function is here \(\sfw(C) \defby \frac{\RCPF_{\Lambda_{\partial}}(x,y)\RCPF_{\Lambda_{\partial}}(\varnothing)\{\mathds{1}_{C_{x,y}(\widehat{\bfn+\bfm})=C}\} }{\RCPF_{\Lambda_{\partial}}(\varnothing)^2}\).

\subsection{\(\beta\geq 0\), \(h>0\)}

This last case is handled in the same way we handled the presence of a boundary condition. Start with \((\Lambda,E_{\Lambda})\) and add a vertex (the \define{ghost}) \(g\) to \(\Lambda\). Then, add an edge between \(i\) and \(g\) for each \(i\in\Lambda\). Denote the obtained augmented graph \((\Lambda_g,E_{\Lambda_g})\). Doing the same expansion as before, one obtains
\begin{equation*}
\sum_{\sigma\in\{-1,1\}^{\Lambda}} \sigma_A  e^{\beta \sum_{ij\in E_{\Lambda}}\sigma_i\sigma_j +h\sum_{i\in\Lambda} \sigma_i } = 2^{|\Lambda|}
\begin{cases}
\RCPF_{\Lambda_{g}}(A) & \text{ if } |A| \text{ even},\\
\RCPF_{\Lambda_{g}}(A\cup\{\partial\}) & \text{ if } |A| \text{ odd},
\end{cases}
\end{equation*}
where \(\RCPF_{\Lambda_{g}}(A) \defby \sum_{\bfn\in\calN_A(\Lambda_g)} \weight_{\beta,h}(\bfn)\) with \(\weight_{\beta,h}(\bfn) \defby \prod_{e\in E_{\Lambda}}\frac{\beta^{\bfn_e}}{\bfn_e!}\prod_{ig\in E_{\Lambda_g}} \frac{h^{\bfn_e}}{\bfn_e!}\).
As a direct consequence,
\begin{equation*}
\langle \sigma_{A}\rangle_{\Lambda,\beta,h} = \frac{1}{\RCPF_{\Lambda_{g}}(\varnothing)}\begin{cases}
\RCPF_{\Lambda_{g}}(A) & \text{ if } |A| \text{ even},\\
\RCPF_{\Lambda_{g}}(A\cup\{g\}) & \text{ if } |A| \text{ odd}.
\end{cases}
\end{equation*}

We can now achieve our goal, that is to obtain a graphical representation of the truncated two-point function. Proceeding as in the case with a boundary (applying the Switching Lemma), one gets
\begin{equation}\label{eq:GeomReprTrunc2PTF}
\langle \sigma_{x};\sigma_y\rangle_{\Lambda} = \frac{\RCPF_{\Lambda_{g}}(x,y)\RCPF_{\Lambda_{g}}(\varnothing)\{\mathds{1}_{y\nleftrightarrow g}\} }{\RCPF_{\Lambda_{g}}(\varnothing)^2} = \sum_{\substack{C\subset \Lambda_g\\ x,y\in C}} \mathds{1}_{g\notin C} \sfw(C) \frac{\RCPF_{\Lambda_{g}}(\varnothing)\{\mathds{1}_{C\nleftrightarrow g}\}}{\RCPF_{\Lambda_{g}}(\varnothing)},
\end{equation}
where \(\sfw(C) \defby \frac{\RCPF_{\Lambda_{g}}(x,y)\{\mathds{1}_{C_{xy}=C }\}}{\RCPF_{\Lambda_{g}}(\varnothing)}\). We thus obtained a graphical representation with the cost of a cluster coming from the fact that we have to disconnect it from the ghost (that models the action of the magnetic field).

\section*{Note added in proofs}
This paper was written already a few years ago and some progress has been made
since then. In particular, concerning the second item in the Open problems section,
the situation is now much better understood. We refer to the following papers for
more information: \cite{Aoun+Ioffe+Ott+Velenik-2021, Aoun+Ott+Velenik-2022, Aoun+Ott+Velenik-2021}.

\bibliographystyle{plain}
\bibliography{OV19-survey}

\end{document}